\documentclass[11pt,oneside,notitlepage]{amsart}
\usepackage{mathpazo}
\usepackage[symbol,perpage]{footmisc}

\vfuzz \hfuzz
\let\Person\textsc

\let\Title\textit

\begin{document}

\begin{center}
\textsc{\Large\bf There are infinitely many prime numbers in
all arithmetic progressions with first term and difference coprime.}\footnote{Originally published in Abhandlungen der K\"{o}niglich Preussischen Akademie der Wissenschaften von 1837, 45--81. Translated by Ralf Stephan,
eMail: \texttt{mailto:ralf@ark.in-berlin.de}. The scanned images are available
at 
\texttt{http://bibliothek.bbaw.de/bibliothek-digital/digitalequellen/schriften/
anzeige?band=07-abh/1837\&seite:int=00000286}}\\[5pt]
(By Mr. \textit{Lejeune-Dirichlet})\\[5pt]
[Read to the Academy of Sciences the 27th of July, 1837]\\[15pt]
\end{center}

\bigskip

MSC-Class: 01A55; 11-03; 11N13; 11B25 

\bigskip

Observant investigation of the series of primes will perceive several
properties with a generality that can be lifted to any degree of
probability by using continued induction, while the discovery of their
proof with all required strictness presents the greatest difficulties.
One of the most curious results of this kind appears when we divide the
members of the series through an arbitrary number. If we ignore the
primes dividing the divisor, which are among the first members of the
series, all other primes will leave a residue that is coprime to the
divisor. The main result after continued division is that every residue
appears to return infinitely often, and in particular, that the ratio of
the values that indicate how often two arbitrary residues have come up
until a specific position is reached will have unity as limit if we continue
division indefinitely. After abstracting away the constancy of appearance
of single residues and limiting the result to the never ending of the
appearance of each residue, we can state the latter as the theorem:
``that each unlimited arithmetic progression, with the first member and the
difference being coprime, will contain infinitely many primes.''

No proof existed for this simple theorem until now, however desirable such
a proof would have been
for the numerous applications that can be made of the
theorem. The only mathematician who has tried a justification of this
theorem, as far as I know, is
\Person{Legendre}\cite{L},
who should not only have been interested in investigating it because
the difficulty of the subject would have appealed to him, but especially
also because he used the mentioned property of arithmetic progressions
as lemma with some earlier works. \Person{Legendre} bases a possible proof
on the task to find the longest run of members of an arithmetic progression
that are divisible by given primes, but he solves it only by induction.
If one tries to prove the solution of that task, which was thus found by him
and which is highly strange because of its simplicity, then one encounters
great difficulties that I did not succeed to overcome. Only after I
entirely left the line taken by \Person{Legendre} I arrived at a completely
strict proof of the theorem on arithmetic progressions. The proof so found,
which I have the honor to present to the Academy in this paper, is not
fully arithmetical but based partly on the study of continous
variables. With the novelty of the applied principles it appeared useful
to me to start with the treatment of the special case where the difference
of the progression is an odd prime, before proving the theorem in its
entire generality.

\bigbreak

\centerline{\S.~1.}

\nobreak\bigskip

Let $p$ an odd prime and $c$ a primitive root of it such that the residues
of the powers:
\[ c^0, c^1, c^2,\ldots, c^{p-2},\]
when divided by $p$, and ordered, are identical to the numbers:
$ 1,2,3,\ldots, p-1$. 
Let $n$ a number not divisible by $p$, then, after \Person{Gauss},
we will call the exponent $\gamma < p - 1$ which satisfies the congruence
$c^\gamma \equiv n$ (mod.~$p$) the index of~$n$ and, if necessary,
denote it as $\gamma_n$. The choice of the primitive root shall be
arbitrary as long as it is fixed. Regarding the above defined indices
the easily provable theorem holds that the index of a product equals
the sum of indices of its factors minus the included multiple of $p-1$.
Further we notice that always 
$\gamma_1 = 0$, $\gamma_{p-1} = \frac{1}{2}(p-1)$, as well as
$\gamma_n$ is even or odd, according to $n$ being quadratic residue of~$p$
or not, or applying \Person{Legendre}'s symbol, according to
$ \bigl( \frac{n}{p} \bigr) = +1$ or
$ \bigl( \frac{n}{p} \bigr) = -1$.

Now let $q$ be a prime different from $p$ (not excluding~$2$) and
$s$ a positive variable greater than unity. Further we denote as~$\omega$
an arbitrary root of the equation:
\begin{equation}
\label{eqn-1}
\omega^{p-1} - 1 = 0,
\end{equation}
and we construct the geometrical series:
\begin{equation}
\label{eqn-2}
\frac{1}{1 - \omega^{\gamma} \frac{1}{q^s}}
   = 1 + \omega^\gamma \frac{1}{q^s}
       + \omega^{2\gamma} \frac{1}{q^{2s}}
       + \omega^{3\gamma} \frac{1}{q^{3s}}
       + \cdots,
\end{equation}
where $\gamma$ means the index of $q$. If we substitute for~$q$
every prime different from~$p$ and multiply the so formed equations
in each other we get a series on the right hand with a structure that
is easily perceived. Namely, let $n$ any integer not divisible by~$p$,
and let $n = q'^{m'}  q''^{m''} \ldots$,
where $q'$, $q'',\ldots$ denote different primes,
then the general term will be of the form:
\[ \omega^{m'\gamma_{q'} + m''\gamma_{q''} + \cdots} \frac{1}{n^s}.\]
But now it holds that:
\[ m' \gamma_{q'} + m'' \gamma_{q''} + \cdots \equiv \gamma_n
   \mbox{ (mod.\ $p - 1$)},\]
and because of (\ref{eqn-1}):
\[ \omega^{m'\gamma_{q'} + m''\gamma_{q''} + \cdots}
   = \omega^{\gamma_n}.\]

Therefore we have the equation:
\begin{equation}
\label{eqn-3}
\prod \frac{1}{1 - \omega^{\gamma} \frac{1}{q^s}}
   = \sum \omega^{\gamma} \frac{1}{n^s} = L,
\end{equation}
where the multiplication sign applies to the whole series of primes
with the only exception of~$p$, while the sum is over all integers
from $1$ to~$\infty$ not divisible by~$p$. The letter~$\gamma$ means
on the left $\gamma_q$, but on the right $\gamma_n$.

The equation just found represents $p - 1$ different equations that
result if we put for $\omega$ its $p - 1$ values. It is known that
these $p - 1$ different values can be written using powers of the
same $\Omega$ when it is chosen correctly, to wit:
\[ \Omega^0, \Omega^1, \Omega^2,\ldots, \Omega^{p-2}.\]

According to this notation, we will write the different values~$L$
of the series or product as:
\begin{equation}
\label{eqn-4}
L_0, L_1, L_2,\ldots, L_{p-2}
\end{equation}
where it is obvious that $L_0$ and
$L_{\frac{1}{2} (p - 1)}$ have a meaning independent of the choice of~$\Omega$
and that they relate to $\omega = 1$ and $\omega = -1$, respectively.

Before we go on it is necessary to state the reason for the condition
made above, that $s > 1$ should hold. We can convince ourselves of the
necessity of this limitation if we respect the essential difference
which exists between two kinds of infinite series. If we regard each
value instead of each term or, it being imaginary, its module, then two
cases can happen. Either it is possible to give a finite value which is
greater than the sum of any of however many of these values or moduli,
or this condition cannot be satisfied by any finite number. In the first
case, the series always converges and has a completely defined sum
regardless how the series terms are ordered, be it that they continue
to two and more dimensions or that they comprise a double or multiple
series. In the second case the series \textbf{can}
converge too but convergence is
essentially dependent on the kind of order of terms. Does convergence
hold for a specific order then it can stop when this order is changed,
or, if this does not happen, then the sum of the series might become
completely different. So, for example, of the two series made from
the same terms:
\begin{eqnarray*}
1 - \frac{1}{\sqrt{2}} + \frac{1}{\sqrt{3}}
  - \frac{1}{\sqrt{4}} + \frac{1}{\sqrt{5}}
  - \frac{1}{\sqrt{6}} + \cdots,\\
1 + \frac{1}{\sqrt{3}} - \frac{1}{\sqrt{2}}
  + \frac{1}{\sqrt{5}} + \frac{1}{\sqrt{7}}
  - \frac{1}{\sqrt{4}} + \cdots,
\end{eqnarray*}
only the first converges while of the following:
\begin{eqnarray*}
1 - \frac{1}{2} + \frac{1}{3}
  - \frac{1}{4} + \frac{1}{5}
  - \frac{1}{6} + \cdots,\\
1 + \frac{1}{3} - \frac{1}{2}
  + \frac{1}{5} + \frac{1}{7}
  - \frac{1}{4} + \cdots,
\end{eqnarray*}
both converge, but with different sums.

Our infinite series~$L$, as can be easily seen, belongs only then to
the first of the classes we just differentiated if we assume
$s > 1$, such that if we put, under this condition,
$L = \lambda + \mu \sqrt{-1}$, then $\lambda$ and $\mu$ will get
completely defined finite values. Denoting now with
$ f_m + g_m \sqrt{-1} $
the product of the first $m$ factors of the form
$ \frac{1}{1 - \omega^{\gamma} \frac{1}{q^s}}$,
where the order of the factors may be arbitrary, we can always make
$m$ of a size such that among these first $m$ factors will be all
those satisfying $q < h$ with $h$ an arbitrary integer. As soon as
$m$ reaches this size, each of both differences $f_m - \lambda$,
$g_m - \mu$, will obviously, ignoring the sign, always stay smaller
than $\frac{1}{h^s} + \frac{1}{(h+1)^s} + \cdots$,
however large $m$ may be imagined to grow further. Under the assumption
$s > 1$ however, the value
$\frac{1}{h^s} + \frac{1}{(h+1)^s} + \cdots $
may shrink arbitrarily small with a correspondingly huge $h$.
Therefore, it is proved that the infinite product in (\ref{eqn-3})
has a value equal to the series $L$, independent of the order of factors.
With $s = 1$ or $s < 1$ however, this proof can be no longer applied
and, in fact, the infinite product has then in general no longer a
definite value, regardless of factor ordering. If we could prove the
existence of a limit for the multiplication continued to infinity,
given a specific ordering of the factors, then the equation (\ref{eqn-3}),
understood correctly, would still hold but would have no use for
the statement of the value. We would then, $q', q'', q'',\ldots$
being the values of~$q$ according to the assumed ordering, have to
view the series~$L$ as a multiple series that has to be ordered such
that first those members would be taken where $n$ only contains the
prime factor~$q'$, then those of the rest where $n$ has no other factors
than~$q', q''$ and so on. From the necessity to order the members this
way the summation of the series would become as difficult as the
investigation of the product itself already is.

\bigbreak

\centerline{\S.~2.}

\nobreak\bigskip

If we put $s = 1 + \varrho$ the equation (\ref{eqn-3}) still holds,
however small the positive value $\varrho$ is assumed. We want to
study now how the series $L$ in (\ref{eqn-3}) changes if $\varrho$
is allowed to become infinitely small. The behaviour of the series with respect
to this is quite different, according to $\omega$ being equal to
positive unity or having any other value. To begin with the first case,
or the investigation of~$L_0$, we look at the sum:
\[ S = \frac{1}{k^{1 + \varrho}}
       + \frac{1}{(k+1)^{1 + \varrho}}
       + \frac{1}{(k+2)^{1 + \varrho}} + \cdots,\]
where $k$ denotes a positive constant. If we substitute in the well-known
formula:
\[ \int_0^1 x^{k-1} \log^\varrho \left( \frac{1}{x} \right) \, dx
   = \frac{\Gamma(1 + \varrho)}{k^{1 + \varrho}} \]
for $k$ sequentially $k,\, k + 1,\, k + 2,\ldots$ and add, we get:
\[ S = \frac{1}{\Gamma(1 + \varrho)} \int_0^1 \log^\varrho \left(
      \frac{1}{x} \right) \frac{x^{k-1}}{1 - x}\, dx.\]
If we add $\frac{1}{\varrho}$ and subtract at the same time:
\[ \frac{1}{\varrho} = \frac{\Gamma(\varrho)}{\Gamma(1 + \varrho)}
   = \frac{1}{\Gamma(1 + \varrho)}
      \int_0^1 \log^{\varrho - 1} \left( \frac{1}{x} \right) \, dx,\]
the equation transforms into:
\[ S = \frac{1}{\varrho} + \frac{1}{\Gamma(1 + \varrho)}
      \int_0^1 \left( \frac{x^{k-1}}{1 - x}
            - \frac{1}{\log \left( \frac{1}{x} \right)}
            \right) \log^\varrho \left( \frac{1}{x} \right) \, dx,\]
where the second member, for $§\varrho$ infinitely small, approaches
the finite limit:
\[\int_0^1 \left( \frac{x^{k-1}}{1 - x}
            - \frac{1}{\log \left( \frac{1}{x} \right)}
            \right) \, dx .\]

Regarding instead of the series $S$ the more general one which has two
positive constants $a$, $b$:
\[ \frac{1}{b^{1 + \varrho}}
   + \frac{1}{(b + a)^{1 + \varrho}}
   + \frac{1}{(b + 2a)^{1 + \varrho}}
   + \cdots,\]
we need only transform it into:
\[ \frac{1}{a^{1 + \varrho}} \Biggl(
   \frac{1}{\bigl( \frac{b}{a} \bigr)^{1 + \varrho}}
   + \frac{1}{\bigl( \frac{b}{a} + 1 \bigr)^{1 + \varrho}}
   + \frac{1}{\bigl( \frac{b}{a} + 2 \bigr)^{1 + \varrho}}
   + \cdots \Biggr) \]
and compare with $S$ to see immediately that it equals to an expression
of the following form:
\[ \frac{1}{a} \cdot \frac{1}{\varrho} + \varphi(\varrho),\]
where $\varphi(\varrho)$ approaches a finite limit with $\varrho$ becoming
infinitely small.

The studied series $L_0$ consists of $p - 1$ partial series like:
\[ \frac{1}{m^{1 + \varrho}}
   + \frac{1}{(p + m)^{1 + \varrho}}
   + \frac{1}{(2p + m)^{1 + \varrho}} + \cdots,\]
where we have to assume successively $m = 1,\,2,\ldots,\, p-1$.
We have thus:
\begin{equation}
\label{eqn-5}
L_0 = \frac{p-1}{p}\cdot\frac{1}{\varrho} + \varphi(\varrho),
\end{equation}
where again $\varphi(\varrho)$ is a function of $\varrho$ that, whenever
$\varrho$ gets infinitely small, has a finite value which could be easily
expressed through a definite integral, given what we found so far. This
is not necessary for our task, however. The equation (\ref{eqn-5}) shows
that, for infinitely small $\varrho$, $L_0$ will become $\infty$ such that
$L_0 - \frac{p-1}{p}\cdot\frac{1}{\varrho}$
remains finite.

\bigbreak

\centerline{\S.~3.}

\nobreak\bigskip

After we have found how our series, with $\omega = 1$ assumed, changes
with $s$ approaching unity from above it remains
for us to extend the same study
to the other roots $\omega$ of the equation $\omega^{p-1} - 1 = 0$.
Although the sum of the series~$L$ is independent of the ordering of its
members, as long as $s > 1$, it will still be of advantage if we imagine
the members ordered such that the values of $n$ will continously increase.
On this condition,
\[ \sum \omega^\gamma \frac{1}{n^s} \]
will be a function of~$s$ that remains continous and finite for all
positive values of~$s$. Thus the limit that is approached by the value
of the series if $s = 1 + \varrho$ and $\varrho$ is let become infinitely
small, and which is independent of the ordering of the members is expressed
by:
\[ \sum \omega^\gamma \frac{1}{n} \]
which wouldn't necessarily be with a different ordering as
$\sum \omega^{\gamma} \frac{1}{n}$ would then differ from
$\sum \omega^{\gamma} \frac{1}{n^{1 + \varrho}}$
by a finite amount or might not even have a value.

To prove the statement just made we denote as $h$ an arbitrary positive
integer and express the sum of the first $h(p-1)$ members of the series:
\[ \sum \omega^\gamma \frac{1}{n^s} \]
with the help of the formula already used above, which holds for any
positive $s$:
\[ \int_0^1 x^{n-1} \log^{s-1} \left( \frac{1}{x} \right) \, dx
   = \frac{\Gamma(s)}{n^s} \]
by a definite integral. We therefore get for the sum:
\[ \frac{1}{\Gamma(s)} \int_0^1
   \frac{\frac{1}{x} \, f(x)}{1 - x^p}
   \log^{s-1} \left( \frac{1}{x} \right) \,dx
 - \frac{1}{\Gamma(s)} \int_0^1
   \frac{\frac{1}{x} \, f(x)}{1 - x^p}
   x^{hp} \log^{s-1} \left( \frac{1}{x} \right) \,dx \]
where we have used the abbreviation:
\[ f(x) = \omega^{\gamma_1} x + \omega^{\gamma_2} x^2 +
      \cdots + \omega^{\gamma_{p-1}} x^{p-1}. \]
If we assume now $\omega$ not $= 1$, then the polynomial
$ \frac{1}{x} \, f(x)$ is divisible by $1 - x$ because we have:
\[ f(1) = \omega^{\gamma_1} + \omega^{\gamma_2} +
      \cdots + \omega^{\gamma_{p-1}}
   = 1 + \omega + \cdots + \omega^{p-2} = 0.\]
If we eliminate thus the factor $1-x$ from numerator and denominator of
the fraction under the integral sign the fraction becomes:
$$\frac{t + u\sqrt{-1}}{1 + x + x^2 + \cdots + x^{p-1}},$$
where $t$ and $u$ are polynomials with real coefficients. If we write
$T$ and $U$ for the largest possible values of $t$ and $u$ between $x = 0$
and $x = 1$ then obviously the real and imaginary parts of the second
integral:
\[ \frac{1}{\Gamma(s)} \int_0^1
   \frac{\frac{1}{x} \, f(x)}{1 - x^p}
   x^{hp} \log^{s-1} \left( \frac{1}{x} \right) \,dx \]
are smaller than
\begin{eqnarray*}
\frac{T}{\Gamma(s)} \int_0^1
      x^{hp} \log^{s-1} \left( \frac{1}{x} \right) \, dx
   &=& \frac{T}{(hp + 1)^s},\\
\frac{U}{\Gamma(s)} \int_0^1
      x^{hp} \log^{s-1} \left( \frac{1}{x} \right) \, dx
   &=& \frac{U}{(hp + 1)^s},\\
\end{eqnarray*}
respectively, and so the integral disappears for $h = \infty$. The series:
\[ \sum \omega^\gamma \frac{1}{n^s}, \]
with the assumed ordering of its members, converges therefore, and for
its sum we have the expression:
\[ \sum \omega^\gamma \frac{1}{n^s}
   = \frac{1}{\Gamma(s)} \int_0^1
      \frac{\frac{1}{x} \, f(x)}{1 - x^p}
      \log^{s-1} \left( \frac{1}{x} \right) \,dx.\]
This function of $s$ not only itself remains continous and finite as long
as~$s > 0$, but the same property applies also to its derivative with respect
to~$s$. To convince oneself of this it is enough to remember that
$\Gamma(s)$, $\frac{d\Gamma(s)}{ds}$, is continous and finite
too, and that $\Gamma(s)$ never disappears as long as $s$ remains positive.

Thus, if we put:
\[\frac{1}{\Gamma(s)} \int_0^1
   \frac{\frac{1}{x} \, f(x)}{1 - x^p}
   \log^{s-1} \left( \frac{1}{x} \right) \,dx
  = \psi(x) + \chi(s)\sqrt{-1},\]
where $\psi(s)$ und $\chi(s)$ are real functions, we have for positive
$\varrho$, after a well known theorem:
\begin{equation}
\label{eqn-6}
\psi(1 + \varrho)
   = \psi(1) + \varrho \psi'(1 + \delta \varrho),\qquad
\chi(1 + \varrho)
   = \chi(1) + \varrho \chi'(1 + \varepsilon \varrho),
\end{equation}
where we abbreviated:
\[ \psi'(s) = \frac{d\psi(s)}{ds},\qquad
   \chi'(s) = \frac{d\chi(s)}{ds}\]
and denoted as $\delta$ and $\varepsilon$ positive fractions independent
of $\varrho$.

Incidentally, it is easily understandable that with $\omega = -1$ we get:
$ \chi(s) = 0$, and that if we go from an imaginary root $\omega$ to its
conjugate $ \frac{1}{\omega}$ then $\psi(s)$ will have the
same value while $\chi(s)$ will become its negative.

\bigbreak

\centerline{\S.~4.}

\nobreak\bigskip

We have to prove now that the finite limit approached by
$ \sum \omega^\gamma \frac{1}{n^{1+\varrho}}$,
with the positive $\varrho$ becoming infinitely small, and given that
$\omega$ does not mean the root~$1$, will be nonzero. This limit is,
after the last section $\sum \omega^\gamma \frac{1}{n}$
and expressed by the integral $\sum \omega^\gamma \frac{1}{n} 
   = - \int_0^1 \frac{\frac{1}{x} \, f(x)}{x^p - 1} \,dx $
which can be easily written using logarithms and circular functions.

Let us take an arbitrary linear factor of the denominator $x^p - 1$:
\[ x - e^{\frac{2m\pi}{p} \sqrt{-1}},\]
where $m$ is of the series $0,\, 1,\, 2,\ldots,\, p-1$. If we decompose:
\[ \frac{\frac{1}{x} \, f(x)}{x^p - 1} \]
into partial fractions then, after known formulae, the numerator of
the fraction:
\[ \frac{A_m}{x - e^{\frac{2m\pi}{p} \sqrt{-1}}} \]
is given by the expression:
\[ \frac{\frac{1}{x} \, f(x)}{px^{p-1}} \]
where $x = e^{\frac{2m\pi}{p} \sqrt{-1}}$.  So we have:
\[ A_m = \frac{1}{p} f \left( e^{\frac{2m\pi}{p} \sqrt{-1}} \right).\]
If we substitute this value and note that $A_0 = 0$ we get:
\[ \sum \omega^\gamma \frac{1}{n}
  = - \frac{1}{p} \sum f \left( e^{\frac{2m\pi}{p} \sqrt{-1}} \right)
      \int_0^1 \frac{dx}{x - e^{\frac{2m\pi}{p} \sqrt{-1}}},\]
where the sum on the right goes from $m = 1$ to $m = p - 1$.

The function:
\[ f \left( e^{\frac{2m\pi}{p} \sqrt{-1}} \right) \]
is well known from cyclotomy and can be easily related to:
\[ f \left( e^{\frac{2\pi}{p} \sqrt{-1}} \right). \]
It namely holds that
\[ f \left( e^{\frac{2m\pi}{p} \sqrt{-1}} \right)
   = \sum \omega^{\gamma_g} e^{gm \frac{2\pi}{p} \sqrt{-1}},\]
where the sum is from $g = 1$ to $g = p - 1$. If we substitute for
$gm$ the respective residue $h$ modulo~$p$ then $1,\, 2,\ldots,\, p - 1$
become the different values of~$h$, and we have, because of
$gm \equiv h$ (mod.~$p$):
\[ \gamma_g \equiv \gamma_h - \gamma_m \mbox{ (mod.\ $p - 1$)}.\]

Thus if we write $\gamma_h - \gamma_m$ for $\gamma_g$, which is allowed
because of the equation $\omega^{p-1} - 1 = 0$, then we get:
\[ f \left( e^{\frac{2m\pi}{p} \sqrt{-1}} \right)
   = \omega^{-\gamma_m} \sum \omega^{\gamma_h}
         e^{h \frac{2\pi}{p} \sqrt{-1}}
   = \omega^{-\gamma_m}
      f \left( e^{\frac{2\pi}{p} \sqrt{-1}} \right).\]
The equation above becomes therefore:
\[ \sum \omega^\gamma \frac{1}{n}
   = - \frac{1}{p} f \left( e^{\frac{2\pi}{p} \sqrt{-1}} \right)
      \sum \omega^{-\gamma_m} 
      \int_0^1 \frac{dx}{x - e^{\frac{2m\pi}{p} \sqrt{-1}}}.\]
Now, for any positive fraction $\alpha$:
\[ \int_0^1 \frac{dx}{x - e^{2\alpha\pi \sqrt{-1}}}
   = \log(2 \sin \alpha \pi)
      + \frac{\pi}{2}(1 - 2\alpha) \sqrt{-1},\]
therefore:
\[ \sum \omega^\gamma \frac{1}{n}
   = - \frac{1}{p} f \left( e^{\frac{2\pi}{p} \sqrt{-1}} \right)
      \sum \omega^{-\gamma_m} \left(
     \log \left( 2 \sin \frac{m\pi}{p} \right)
      + \frac{\pi}{2} \left( 1 - \frac{2m}{p} \right) \sqrt{-1}
      \right).\]

Although this expression is very simple for
$\sum \omega^\gamma \frac{1}{n}$,
in general we cannot conclude that
$\sum \omega^\gamma \frac{1}{n}$
has a nonzero value. What is missing are fitting principles for
the statement of conditions under which transcendent compounds containing
undefined integers can disappear. But our desired proof succeeds for
the specific case where $\omega = -1$. For imaginary values of $\omega$
we will give another method in the following section that, however, cannot
be applied to the mentioned specific case.
On the condition $\omega = -1$, and with $\gamma_m$ even or odd according to:
$  \bigl( \frac{m}{p} \bigr) = +1$ or $-1$,
and thus $ (-1)^{-\gamma_m} = \bigl( \frac{m}{p} \bigr) $,
as well as $ (-1)^{\gamma_n} = \bigl( \frac{n}{p} \bigr)$,
we get as limit of $L_{\frac{1}{2}(p-1)}$ for $\varrho$ becoming infinitely
small:
\begin{multline*}
\sum \left( \frac{n}{p} \right) \frac{1}{n} \\
   = - \frac{1}{p} f \left( e^{\frac{2\pi}{p} \sqrt{-1}} \right)
      \sum \left( \frac{m}{p} \right) \left(
     \log \left( 2 \sin \frac{m\pi}{p} \right)
      + \frac{\pi}{2} \left( 1 - \frac{2m}{p} \right) \sqrt{-1}
      \right),
\end{multline*}
or more simple, since
$\sum \bigl( \frac{m}{p} \bigr) = 0$
if we sum from $m = 1$ to $m = p - 1$:
\[ \sum \left( \frac{n}{p} \right) \frac{1}{n}
   = - \frac{1}{p} f \left( e^{\frac{2\pi}{p} \sqrt{-1}} \right)
      \sum \left( \frac{m}{p} \right) \left(
     \log \left( 2 \sin \frac{m\pi}{p} \right)
      - \frac{\pi}{p} m \sqrt{-1} \right).\]

We have to distinguish two cases, depending on the prime $p$ having
the form $4\mu + 3$ or $4\mu + 1$.  In the first case it holds for two
values, like $m$ und $p - m$, adding to $p$ that:
\[ \left( \frac{m}{p} \right)
   = - \left( \frac{p-m}{p} \right)
      \mbox{ and }
\sin \frac{m\pi}{p} = \sin \frac{(p-m)\pi}{p}.\]
Therefore the real part of the sum disappears and we get, denoting with
$a$ those values of $m$ for which $ \bigl( \frac{m}{p} \bigr) = 1$
and with $b$ those for which $ \bigl( \frac{m}{p} \bigr) = -1$,
or in other words, $a$ and $b$ denoting the quadratic residues and
nonresidues of $p$ that are smaller than $p$:
\[ \sum \left( \frac{n}{p} \right) \frac{1}{n}
   = \frac{\pi}{p^2} f \left( e^{\frac{2\pi}{p} \sqrt{-1}} \right)
      \bigl(\sum a - \sum b\bigr) \sqrt{-1}.\]
With $p = 4\mu + 1$ the imaginary part of the sum disappears because then
$ \bigl( \frac{m}{p} \bigr) = \bigl( \frac{p-m}{p} \bigr)$,
and we get:
\[ \sum \left( \frac{n}{p} \right) \frac{1}{n}
   = \frac{1}{p} f \left( e^{\frac{2\pi}{p} \sqrt{-1}} \right)
      \log \frac{\prod \sin \frac{b\pi}{p}}{
                 \prod \sin \frac{a\pi}{p}},\]
where the multiplication extends over all $a$ or $b$.

Notice know that, under the assumption $\omega = -1$, using well-known
formulae\cite{DD}
we get for $f \left( e^{\frac{2\pi}{p} \sqrt{-1}} \right)$ in the first
case $\sqrt{p}\sqrt{-1}$, and in the last $\sqrt{p}$, so respectively:
\[ \sum \left( \frac{n}{p} \right) \frac{1}{n}
      = \frac{\pi}{p\sqrt{p}} (\sum b - \sum a),\qquad
   \sum \left( \frac{n}{p} \right) \frac{1}{n}
      = \frac{1}{\sqrt{p}}
         \log \frac{\prod \sin \frac{b\pi}{p}}{
                    \prod \sin \frac{a\pi}{p}}.\]

In the case of $p = 4\mu + 3$ we see immediately that
$ \sum \bigl( \frac{n}{p} \bigr) \frac{1}{n}$
is nonzero since
$\sum a + \sum b = \frac{1}{2} p(p-1)$
is odd and thus $\sum a = \sum b$ cannot hold. To prove the same for
$p = 4\mu + 1$ we use the equations which are known from
cyclotomy\cite{G}
\[ 2 \prod \left( x - e^{\frac{2a\pi}{p}\sqrt{-1}} \right)
      = Y - Z\sqrt{p},\qquad
   2 \prod \left( x - e^{\frac{2b\pi}{p}\sqrt{-1}} \right)
      = Y + Z\sqrt{p},\]
where $Y$ and $Z$ are polynomials with integer coefficients. If we
substitute in these equations and those following from them:
\[ 4 \frac{x^p - 1}{x - 1} = Y^2 - p Z^2 \]
$x = 1$ and denote then with $g$ and $h$ those integer values taken by $Y$
and $Z$, we get after several easy reductions:
\[ 2^{\frac{p+1}{2}} \prod \sin \frac{a\pi}{p}
     = g - h\sqrt{p},\quad
   2^{\frac{p+1}{2}} \prod \sin \frac{b\pi}{p}
     = g + h\sqrt{p},\quad
   g^2 - p h^2 = 4p.\]
From the last equation follows that $g$ is divisible by $p$. If we
therefore put $g = pk$ and divide the first two [equations] through
each other we get:
\[ \frac{\prod \sin \frac{b\pi}{p}}{
         \prod \sin \frac{a\pi}{p}}
      = \frac{k\sqrt{p} + h}{k\sqrt{p} - h},\qquad
   h^2 - p k^2 = -4.\]

From the second of these equations, $h$ cannot be zero, thus both sides
of the first [equation] are different from unity, from which immediately
follows, respecting the expression we obtained above, that:
$ \sum \bigl( \frac{n}{p} \bigr) \frac{1}{n}$
cannot have the value zero, q.~e.~d.

We can add that the sum
$ \sum \bigl( \frac{n}{p} \bigr) \frac{1}{n}$,
because as limit of a product of only positive factors, namely as limit of:
\[ \prod \frac{1}{ 1 - \left( \frac{q}{p} \right)
      \frac{1}{q^{1 + \varrho}}},\]
for $\varrho$ becoming infinitely small, it can never become negative,
so it will be necessarily positive.

From this statement two important theorems follow directly that are
probably not provable otherwise, of which the one related to the case
$p = 4\mu + 3$ is that for a prime of this form $\sum b > \sum a$ always
holds. However, we don't want to stay with these results of our method
since we will have occasion to get back to the subject with another
investigation.

\bigbreak

\centerline{\S.~5.}

\nobreak\bigskip

To prove for $L_m$, if $m$ is neither $0$ nor $\frac{1}{2}(p-1)$, that
its limit, which corresponds to $\varrho$ being infinitely small, is
different from zero we take the logarithm of:
\[ \prod \frac{1}{ 1 - \omega^\gamma
      \frac{1}{q^{1 + \varrho}}} \]
and develop the logarithm of each factor using the formula:
\[ - \log (1 - x) = x + \frac{1}{2} x^2 + \frac{1}{3} x^3 + \cdots .\]
We so find:
\[ \sum \omega^\gamma \frac{1}{q^{1 + \varrho}}
   + {\textstyle\frac{1}{2}} \sum \omega^{2\gamma}
      \frac{1}{(q^2)^{1 + \varrho}}
   + {\textstyle\frac{1}{3}} \sum \omega^{3\gamma}
      \frac{1}{(q^3)^{1 + \varrho}}
   + \cdots = \log L,\]
where the summation is with repect to $q$ and $\gamma$ means the index
of~$q$. If we substitute for $\omega$ its values:
\[ 1, \Omega, \Omega^2,\ldots, \Omega^{p-2}, \]
add and remember that the sum:
\[ 1 + \Omega^{h\gamma} + \Omega^{2h\gamma} +
   \cdots + \Omega^{(p-2)h\gamma} \]
always disappears except when $h\gamma$ is divisible by $p - 1$
but has in this case the value $p - 1$, and that the condition
$h\gamma \equiv 0$ (mod.\ $p - 1$) is identical with $q^h \equiv 1$
(mod.~$p$), then we get:
\[ (p - 1) \left(
   \sum \frac{1}{q^{1 + \varrho}}
   + {\textstyle\frac{1}{2}} \sum \frac{1}{q^{2 + 2\varrho}}
   + {\textstyle\frac{1}{3}} \sum \frac{1}{q^{3 + 3\varrho}}
   + \cdots \right) = \log (L_0 L_1 \ldots L_{p-2}),\]
where the first, second,\dots\ summation relates to those values of~$q$,
the first, second,\dots\ powers of which are contained in the form
$\mu p +1$, respectively. Because the left side is real it follows that
the product under the log sign is positive, which is also obvious
otherwise, and that we have to take for the logarithm its arithmetical
value that has no ambiguity. The series on the left hand remains always
positive, and we will show now that the right side to the contrary would
be $-\infty$ with vanishing $\varrho$, if we would try to assume the limit for
$L_m$ to disappear.
The right side can be written as:
\[ \log L_0 + \log L_{\frac{1}{2}(p-1)}
   + \log L_1 L_{p-2} + \log L_2 L_{p-3} + \cdots,\]
where $\log L_0$ after \eqref{eqn-5} is equal to the expression:
\[ \log \left( \frac{p-1}{p}\cdot\frac{1}{\varrho} + \varphi(\varrho)
      \right) \]
or:
\[ \log \left( \frac{1}{\varrho} \right) + \log \left(
   \frac{p-1}{p} + \varrho \varphi(\varrho) \right), \]
the second term of which approaches the finite limit:
$ \log \bigl( \frac{p-1}{p} \bigr);$
likewise, $\log L_{\frac{1}{2}(p-1)}$ remains finite, since the limit
of $L_{\frac{1}{2}(p-1)}$ is nonzero with \S.~4. One of the other logarithms
$\log L_m L_{p-1-m}$, is after \S.~3:
\[ \log \bigl(\psi^2(1 + \varrho) + \chi^2 (1 + \varrho)\bigr),\]
which expression, if $L_m$ and thus $L_{p-1-m}$ too would have zero as
limit such that at the same time $\psi(1) = 0$, $\chi(1) = 0$,
would transform into:
\[ \log \bigl(\varrho^2 (\psi'^2 (1 + \delta \varrho)
      + \chi'^2 (1 + \varepsilon \varrho)\bigr) \]
or:
\[ -2 \log \left( \frac{1}{\varrho} \right)
      + \log \bigl(\psi'^2 (1 + \delta \varrho)
      + \chi'^2 (1 + \varepsilon \varrho)\bigr). \]
Combining the term
$ - 2 \log \bigl( \frac{1}{\varrho} \bigr)$
with the first term of $\log L_0$, there results:
$ - \log \bigl( \frac{1}{\varrho} \bigr)$,
which value would become $-\infty$ with infinitely small $\varrho$,
and it is clear that this infinitely large negative value cannot be
cancelled by e.g.:
\[ \log \bigl(\psi'^2(1 + \delta \varrho)
      + \chi'^2 (1 + \varepsilon \varrho)\bigr) \]
because this expression either remains finite or itself becomes
$- \infty$, namely when simultaneously $\psi'(1) = 0$, $\chi'(1) = 0$.
Just as evident is that, would we try to view some other pairs of
related $L$ than $L_m$ and $L_{p-1-m}$ as mutually cancelling, the
contradiction would be even intensified. Therefore it is proved that
the limit of $L_m$ for $m>0$, corresponding to infinitely small
$\varrho$, is finite and different from zero. Also, $L_0$ in the
same case, becomes $\infty$ from which immediately follows that the
series:
\begin{equation}
\label{eqn-7}
\sum \omega^\gamma \frac{1}{q^{1 + \varrho}}
   + {\textstyle\frac{1}{2}} \sum \omega^{2\gamma}
      \frac{1}{q^{2 + 2\varrho}}
   + {\textstyle\frac{1}{3}} \sum \omega^{3\gamma}
      \frac{1}{q^{3 + 3\varrho}}
   + \cdots = \log L
\end{equation}
always approaches a finite limit if $\omega$ not $ = 1$, but becomes
infinitely large for $\omega = 1$ if we let $\varrho$ become infinitely small.

Would we want to have the limit itself, which is not necessary for our
task however, its calculation (for $\omega$ not $-1$) using the expression
$\log \bigl(\psi(1) + \chi(1) \sqrt{-1}\bigr)$
would be afflicted with an ambiguity
that could be lifted easily with specialisation, i.e., as soon as
$p$ and $\omega$ is given numerically. If we equal the series \eqref{eqn-7} 
with $u + v \sqrt{-1}$ and therefore:
\[ u + v\sqrt{-1} = \log L
   = \log \bigl(\psi(1 + \varrho) + \chi(1 + \varrho)\sqrt{-1} \bigr),\]
we have:
\[ u = {\textstyle\frac{1}{2}} \log \bigl(\psi^2(1 + \varrho)
         + \chi^2 (1 + \varrho)\bigr),\]
\[ \cos v = \frac{\psi(1 + \varrho)}{\sqrt{
   \psi^2 (1 + \varrho) + \chi^2(1 + \varrho)}},\qquad
   \sin v = \frac{\chi(1 + \varrho)}{\sqrt{
   \psi^2 (1 + \varrho) + \chi^2(1 + \varrho)}},\]
and thus the limit of $u$ is no longer ambigous:
\[ {\tfrac{1}{2}} \log \bigl(\psi^2(1) + \chi^2(1)\bigr).\]
To get the same with $v$ we note that the series, however small $\varrho$
may be, is continously variable with respect to this value, which can be
easily proved. Therefore also $v$ is a continous function of~$\varrho$.
Because $\psi(1) = 0$, $\chi(1) = 0$ cannot hold at the same time, it will
be possible to derive from the expressions of $\psi(1 + \varrho)$ and
$\chi(1 + \varrho)$ given above as definite integrals always a positive
finite value~$R$ of such character that at least one of the functions
$\psi(1 + \varrho)$, $\chi(1 + \varrho)$ retains its sign for each $\varrho$
smaller than $R$. Therefore $\cos v$ or $\sin v$, as soon as $\varrho$
decreases below $R$, will no longer change signs and thus the continously
variable arc $v$ will no longer be able to increase or decreases by $\pi$.
If we thus determine the finite value of~$v$ corresponding to $\varrho = R$,
let's call it~$V$, which we can easily find by numerical computation from
the series~(\ref{eqn-7}) because [the series] for each finite value of
$\varrho$ belongs to one of the classes differentiated in \S.~1 and thus
has a completely defined sum, the limit $v_0$ of $v$ is given by the
equations
\[ \cos v_0 = \frac{\psi(1)}{\sqrt{\psi^2 (1) + \chi^2(1)}},\qquad
   \sin v_0 = \frac{\chi(1)}{\sqrt{\psi^2 (1) + \chi^2(1)}},\]
under the condition that the difference $V - v_0$, ignoring signs,
has to be smaller than $\pi$.

\bigbreak

\centerline{\S.~6.}

\nobreak\bigskip

We are now able to prove that each arithmetic progression with difference
$p$ whose first member is not divisible by $p$ contains infinitely many primes;
or, in other words, that there are infinitely many primes of the form
$\mu p + m$, where $\mu$ is an arbitrary number and $m$ one of the numbers
$1,\, 2,\, 3,\ldots,\, p - 1$. If we multiply the equations contained in
\eqref{eqn-7} that correspond consecutively to the roots:
\[ 1, \Omega, \Omega^2,\ldots, \Omega^{p-2} \]
with:
\[ 1, \Omega^{-\gamma_m}, \Omega^{-2\gamma_m},\ldots,
   \Omega^{-(p-2) \gamma_m} \]
and add we get on the left side:
\begin{eqnarray*}
\qquad & &\sum ( 1 + \Omega^{\gamma - \gamma_m}
                   + \Omega^{2(\gamma - \gamma_m)} + \cdots
                   + \Omega^{(p-2)(\gamma - \gamma_m)}
               ) \frac{1}{q^{1 + \varrho}} \\
       &+&{\textstyle\frac{1}{2}}
          \sum ( 1 + \Omega^{2\gamma - \gamma_m}
                   + \Omega^{2(2\gamma - \gamma_m)} + \cdots
                   + \Omega^{(p-2)(2\gamma - \gamma_m)}
               ) \frac{1}{q^{2 + 2\varrho}} \\
       &+&{\textstyle\frac{1}{3}}
          \sum ( 1 + \Omega^{3\gamma - \gamma_m}
                   + \Omega^{2(3\gamma - \gamma_m)} + \cdots
                   + \Omega^{(p-2)(3\gamma - \gamma_m)}
               ) \frac{1}{q^{3 + 3\varrho}} \\
       &+& \cdots,
\end{eqnarray*}
where the summation is over $q$ and $\gamma$ denotes the index of
$q$. But now it holds that:
\[ 1 + \Omega^{h\gamma - \gamma_m}
     + \Omega^{2(h\gamma - \gamma_m)} + \cdots
     + \Omega^{(p-2)(h\gamma - \gamma_m)}
   = 0,\]
except when $h\gamma - \gamma_m \equiv 0$ (mod.\ $p - 1$), in which case
the sum equals $p - 1$. This congruence however is identical with
$q^h \equiv m$ (mod.~$p$).  We therefore have the equation:
\begin{multline*}
\sum \frac{1}{q^{1 + \varrho}}
      + {\textstyle\frac{1}{2}} \sum \frac{1}{q^{2 + 2\varrho}}
      + {\textstyle\frac{1}{3}} \sum \frac{1}{q^{3 + 3\varrho}}
      + \cdots\\
   = \tfrac{1}{p-1} (\log L_0
      + \Omega^{-\gamma_m} \log L_1
      + \Omega^{-2\gamma_m} \log L_2
      + \cdots
      + \Omega^{-(p-2)\gamma_m} \log L_{p-2}),
\end{multline*}
where the first summation is over all primes~$q$ of form
$\mu p + m$, the second over all primes~$q$ with squares of that form,
the third over all primes~$q$ with cubes of that form etc.
If we assume now $\varrho$ becoming infinitely small, the right side
will become infinitely large through the term $\log L_0$. Thus also the
left hand has to become infinity. But on this side the sum of all terms,
except the first, remains finite because, as is well-known,
\[ {\textstyle\frac{1}{2}} \sum \frac{1}{q^2}
      + {\textstyle\frac{1}{3}} \sum \frac{1}{q^3}
      + \cdots \]
is even finite if we substitute for $q$ not certain primes, as we did, but
all integers larger than~$1$. Thus the series
\[ \sum \frac{1}{q^{1 + \varrho}} \]
has to grow beyond any positive limit; it has to have infinitely many
terms, that means, there are infinitely many primes~$q$ of form
$\mu p + m$, q.~e.~d.

\bigbreak

\centerline{\S.~7.}

\nobreak\bigskip

To extend the previously given proof to arithmetic progressions with
the difference being any composite integer, several theorems from the
theory of residues will be necessary which we want to collect now to
be able to refer to them more easily in the following. Justification
of these results can be looked up in \emph{Disq.\ arith.\ sect.\ III.}
where the subject is treated in depth.

I. The existence of primitive roots is not limited to odd primes~$p$
but also applies to any power §$p^\pi$ of it. With $c$ a primitive
root for the modulus $p^\pi$, the residues of the powers:
\[ c^0, c^1, c^2,\ldots, c^{(p-1)p^{\pi - 1}-1} \]
with respect to [that modulus] are all mutually different and identical
to the series of the numbers smaller than $p^\pi$ and coprime to $p^\pi$.
If we have now an arbitrary number~$n$ not divisible by~$p$, then the
exponent $\gamma_n < (p-1) p^{\pi - 1}$ which satisfies the congruence
\[ c^{\gamma_n} \equiv n \hbox{ (mod.~$p^\pi$)} \] is completely
determined and we shall call it the index of~$n$. Regarding such indices
the easy to prove theorems hold that the index of a product equals the
sum of factor indices minus the largest contained multiple of
$(p - 1)p^{\pi - 1}$, and that $\gamma_n$ is even or odd corresponding to
$ \bigl( \frac{n}{p} \bigr) = +1$ or $-1$.

II. The prime number $2$ shows an essentially different behaviour in the
theory of primitive roots than the odd primes, and we note about this
prime the following, if we ignore the first power of $2$ which is not
important here.

1) For the modulus $2^2$ we have the primitive root $-1$. If we denote
the index of an arbitrary odd number~$n$ with $\alpha_n$ such that then:
\[ (-1)^{\alpha_n} \equiv n \hbox{ (mod.~$4$)},\]
we have $\alpha_n = 0$ oder $\alpha_n = 1$, according to $n$ having the form
$4\mu + 1$ or $4\mu + 3$, and we get the index of a product if we
subtract its largest contained even number from the sum of factor indices.

2) For a modulus of the form $2^\lambda$, with $\lambda \geq 3$, there
does not exist any primitive root anymore, i.e., there is no number
such that the period of its power residues after the divisor~$2^\lambda$
contains every odd number smaller than~$2^\lambda$. It is possible only
to express one half of these numbers as such residues. If we choose any
number of form $8\mu + 5$ or especially~$5$ as basis then the residues
of the powers
\[ 5^0, 5^1, 5^2,\ldots, 5^{2^{\lambda - 2} - 1} \]
modulo $2^\lambda$ are mutually different and coincide with the numbers 
of form $4\mu + 1$ and being smaller than $2^\lambda$. If we therefore
have a number~$n$ of form~$4\mu + 1$ then the congruence:
\[ 5^{\beta_n} \equiv n \hbox{ (mod.~$2^\lambda$)} \]
can always be satisfied by exactly one exponent or index $\beta_n$ that
should be smaller than $2^{\lambda - 2}$. If $n$ is of form $4 \mu + 3$
then this congruence is impossible. Because however on this condition
$-n$ is of form $4\mu + 1$ then in general we will denote as the
index of an odd number~$n$ the completely defined exponent $\beta_n$
that is smaller than $2^{\lambda - 2}$ and satisfies the congruence:
\[ 5^{\beta_n} \equiv \pm n \hbox{ (mod.~$2^\lambda$)} \]
where the upper or lower sign is applied corresponding to $n$ being of
form $4\mu + 1$ or $4\mu + 3$. Because of this double sign the residue
of~$n$ modulo $2^\lambda$ is thus no longer completely determined by
the index $\beta_n$, since the same index matches two residues that
complement to $2^\lambda$. For the so defined indices obviously the
theorems hold that the index of a product equals the sum of factor
indices minus the largest multiple of $2^{\lambda - 2}$ contained
therein, as well as that $\beta_n$ is even or odd according to $n$
being of form $8\mu \pm 1$ or $8\mu \pm 5$. To remove the aforementioned
ambiguity it is sufficient to study not only the index $\beta_n$ related to the modulus~$2^\lambda$ and base~$5$ but also the index $\alpha_n$ related to
the modulus~$4$ and base~$-1$. In that case, according to $\alpha_n = 0$
or $\alpha_n = 1$, we will have to apply the upper or lower sign in:
\[  5^{\beta_n} \equiv \pm n \hbox{ (mod.~$2^\lambda$)}. \]
We could join, if we want, both indices into one formula and write:
\[ (-1)^{\alpha^n} 5^{\beta_n} \equiv n \hbox{ (mod.~$2^\lambda$)}, \]
by which congruence the residue of $n$ modulo $2^\lambda$ will be
completely determined.

III.  Let now:
\[ k = 2^\lambda p^\pi p'^{\pi'} \ldots,\]
where, as in II.~2, $\lambda \geq 3$  and $p, p',\ldots$ denote mutually
different odd primes. If we take any number~$n$ not divisible by the
primes $2, p, p',\ldots$ and the indices:
\[ \alpha_n, \beta_n, \gamma_n, \gamma'_n,\ldots,\]
corresponding to the moduli:
\[ 4, 2^\lambda, p^\pi, p'^{\pi'},\ldots \]
and their primitive roots:
\[ -1, 5, c, c',\ldots \]
then we have the congruences:
\[ (-1)^{\alpha_n}   \equiv n     \hbox{ (mod.~$4$)},\qquad  
      5^{\beta_n}    \equiv \pm n \hbox{ (mod.~$2^\lambda$)},\]
\[    c^{\gamma_n}   \equiv n     \hbox{ (mod.~$p^\pi$)},\qquad
      c'^{\gamma'_n} \equiv n     \hbox{ (mod.~$p^\pi$)},\ldots,\]
by which the residue of $n$ modulo~$k$ is completely determined, which
follows at once from well-known theorems if we
remember that the double sign in the second of these congruences is fixed
by the first. We will call the indices $\alpha_n, \beta_n, \gamma_n,
\gamma'_n,\ldots$ or $\alpha, \beta, \gamma, \gamma',\ldots$
the system of indices for the number.  Because the Indices:
\[ \alpha, \beta, \gamma, \gamma',\ldots \]
or:
\[ 2,\, 2^{\lambda - 2},\, (p-1)p^{\pi - 1},\,
      (p'-1)p'^{\pi' - 1},\ldots, \]
respectively, can have different values it holds that:
\begin{multline}
\label{eqn-8}
2\cdot2^{\lambda-2} (p-1) p^{\pi - 1} (p' - 1) p'^{\pi' - 1} \cdots
   = k \left( 1 - \frac{1}{2} \right)
       \left( 1 - \frac{1}{p} \right)
       \left( 1 - \frac{1}{p'} \right)
       \cdots\\
   = K
\end{multline}
is the number of possible systems of this kind, which agrees with the
well-known theorem that by $K$ is expressed the number of those numbers
smaller than $k$ and coprime to $k$.

\bigbreak

\centerline{\S.~8.}

\nobreak\bigskip

When trying to prove the theorem on arithmetic progressions in its full
generality we note that we can, without loss of it, assume the
difference~$k$ of the progressions as divisible by~$8$ and thus
of the form included with \S.~7, III. Is the theorem proven under this
condition it will obviously hold as well with the difference odd or only
divisible by~$2$ or~$4$. Let $\theta, \varphi,\omega, \omega',\ldots$
be any roots of the equations:
\[ \theta^2 - 1 = 0,\quad
\varphi^{2^{\lambda - 2}} - 1 = 0,\]
\begin{equation}
\label{eqn-9}
\omega^{(p-1) p^{\pi - 1}} - 1 = 0,\quad
\omega'^{(p'-1) p^{\pi' - 1}} - 1 = 0,\ldots,
\end{equation}
and let $q$ be an arbitrary prime not equal $2, p, p',\ldots$.
If we write the equation:
\begin{multline*}
\frac{1}{1
       - \theta^\alpha \varphi^\beta
                \omega^\gamma \omega'^{\gamma'}
                \cdots \frac{1}{q^s}}\\
   = 1 + \theta^\alpha \varphi^\beta
                \omega^\gamma \omega'^{\gamma'}
                \cdots \frac{1}{q^s}
       + \theta^{2\alpha} \varphi^{2\beta}
                \omega^{2\gamma} \omega'^{2\gamma'}
                \cdots \frac{1}{q^{2s}}
       + \cdots,
\end{multline*}
where $s > 1$, and the system of indices
$\alpha, \beta, \gamma, \gamma',\ldots$ is with respect to $q$, and
if we multiply all equations of this form, which we get if we substitute
for $q$ every prime different from $2, p,p',\ldots$, with each other then
we get, remembering the abovementioned properties of indices and
equations~(\ref{eqn-9}):
\begin{equation}
\label{eqn-10}
\prod \frac{1}{1
        - \theta^\alpha \varphi^\beta
                \omega^\gamma \omega'^{\gamma'}
                \cdots \frac{1}{q^s}}
   = \sum \theta^\alpha \varphi^\beta
                \omega^\gamma \omega'^{\gamma'}
                \cdots \frac{1}{n^s}
   = L,
\end{equation}
where the product is over the primes except $2, p,p'\ldots$, and the
sum is over all positive integers not divisible by the primes
$2, p,p'\ldots$. The system of indices
$\alpha, \beta, \gamma, \gamma',\ldots$
is on the left side with respect to~$q$, on the right side to~$n$.
The general equation \eqref{eqn-10}, in which the different roots
$\theta, \varphi, \omega, \omega',\ldots$ can be mutually combined
arbitrarily, apparently contains a number~$K$ of special equations.
To denote the series~$L$ corresponding to each of the combinations in
a comfortable way we can think of the roots of each of these
equations~(\ref{eqn-9}) expressed as powers of one of them. Let
$\Theta = -1, \Phi, \Omega, \Omega',\ldots$
be roots suitable for that purpose, then:
\[ \theta  = \Theta^{\mathfrak{a}},\qquad
   \varphi = \Phi^{\mathfrak{b}},\qquad
   \omega  = \Omega^{\mathfrak{c}},\qquad
   \omega' = \Omega'^{\mathfrak{c}'},\ldots,\]
where:
\[ \mathfrak{a}  < 2,\qquad
   \mathfrak{b}  < 2^{\lambda - 2},\qquad
   \mathfrak{c}  < (p-1) p^{\pi - 1},\qquad
   \mathfrak{c}' < (p'-1) p'^{\pi' - 1},\ldots\]
and, using this notation, denote the series~$L$ with:
\begin{equation}
\label{eqn-11}
L_{\mathfrak{a}, \mathfrak{b}, \mathfrak{c}, \mathfrak{c}' \cdots}.
\end{equation}
The necessity for the condition $s > 1$ in the equation \eqref{eqn-10}
is the same as already developed in~\S.~1.

\bigbreak

\centerline{\S.~9.}

\nobreak\bigskip

The series denoted with 
$L_{\mathfrak{a}, \mathfrak{b}, \mathfrak{c}, \mathfrak{c}' \cdots}$,
of which the number equals $K$, can be divided into the following three
classes according to the different combinations $\theta, \varphi, \omega, 
\omega',\ldots$ of their roots. The first class contains only one series,
namely $L_{0,0,0,0,\ldots}$, that is, the one where:
\[ \theta  = 1,\quad
   \varphi = 1,\quad
   \omega  = 1,\quad
   \omega' = 1, \ldots \]
holds. The second class shall cover all other series with only real
solutions to the equations \eqref{eqn-9} such that therefore to express
those series we have to combine the signs in:
\[ \theta  = \pm 1,\quad
   \varphi = \pm 1,\quad
   \omega  = \pm 1,\quad
   \omega' = \pm 1, \ldots \]
in every possible way excepting only the combination corresponding to the
first class. The third class finally includes all series~$L$ where at least
one of the roots $\varphi, \omega, \omega',\ldots$ is imaginary, and it is
evident that the series of this class come in pairs since the two root
combinations:
\[ \theta, \varphi, \omega, \omega',\ldots;\quad
\frac{1}{\theta}, \frac{1}{\varphi},
\frac{1}{\omega}, \frac{1}{\omega'},\ldots \]
are mutually different given the just mentioned condition. We will now
have to study the behaviour of those series on substitution of
$s = 1 + \varrho$ and letting become $\varrho$ infinitely small.
Let us look first at the series constituting the first class, and clearly,
we can see it as sum of $K$ partial series of form:
\[ \frac{1}{m^{1 + \varrho}}
   + \frac{1}{(k + m)^{1 + \varrho}}
   + \frac{1}{(2k + m)^{1 + \varrho}}
   + \cdots,\]
where $m < k$ and coprime to $k$. Thus the series of this class equals
after~\S.~2:
\begin{equation}
\label{eqn-12}
\frac{K}{k} \cdot \frac{1}{\varrho} + \varphi(\varrho),
\end{equation}
where $\varphi(\varrho)$ remains finite for infinitely small $\varrho$.

Regarding the series of second or third class we find, if we order their
terms such that the values of $n$ are increasing and with $s > 0$, for
them the equation:
\begin{equation}
\label{eqn-13}
\sum \theta^\alpha \varphi^\beta \omega^\gamma \omega'^{\gamma'}
                \cdots \frac{1}{n^s}
   = \frac{1}{\Gamma(s)} \int_0^1
      \frac{\sum \theta^\alpha \varphi^\beta
                \omega^\gamma \omega'^{\gamma'} \ldots
                x^{n-1}}{1 - x^k}
      \log^{s-1} \left( \frac{1}{x} \right) \,dx,
\end{equation}
where the sum on the right hand is over all positive integers $n$ smaller 
than and coprime to~$k$, and $\alpha, \beta, \gamma, \gamma',\ldots$
stands for the system of indices for~$n$. We easily prove that the right
side has a finite value, because for it we only need to mention that the
polynomial
$\sum \theta^\alpha \varphi^\beta
      \omega^\gamma \omega'^{\gamma'} \ldots x^{n-1}$
involves the factor $1 - x$, which illustrates immediately that if we
put $x = 1$ by which the polynomial is transformed to the product:
\begin{multline*}
(1 + \theta)(1 + \varphi + \cdots + \varphi^{2^{\lambda - 2} - 1})
      (1 + \omega + \cdots + \omega^{(p-1)p^{\pi - 1} - 1})\times\\
      \times
      (1 + \omega' + \cdots + \omega'^{(p'-1)p'^{\pi' - 1} - 1})
      \cdots ,
\end{multline*}
at least one of its factors will disappear, which would be impossible
with the root combination corresponding to the first class:
\[ \theta  = 1,\quad
   \varphi = 1,\quad
   \omega  = 1,\quad
   \omega' = 1, \ldots.\]
We convince ourselves as easily that the right side of the equation
\eqref{eqn-13}, as well as its differential quotient with respect
to~$s$, are continous functions of~$s$. It follows at once that, for
$\varrho$ getting infinitely small, each series of second or third class
approaches a finite limit expressed by:
\begin{equation}
\label{eqn-14}
\sum \theta^\alpha \varphi^\beta \omega^\gamma \omega'^{\gamma'}
                \cdots \frac{1}{n}
   = \int_0^1 \frac{\sum \theta^\alpha \varphi^\beta
                \omega^\gamma \omega'^{\gamma'} \ldots
                x^{n-1}}{1 - x^k} \,dx.
\end{equation}
It remains only to prove that this limit is always nonzero.

\bigbreak

\centerline{\S.~10.}

\nobreak\bigskip

Although the limit for an $L$ of the second or third class is easily
expressed using logarithms and circular functions, like in \S.~4,
such an expression would have no use for the desired proof, also not
if $L$ belongs to the second class, even though this case elsewise
is mostly analogous to that studied in the last half of~\S.~4. Let us
just now assume the mentioned property was proved for each $L$ of the
second class. We will now show how the same requirement can be satisfied
for an $L$ of the third class. To this end we take logarithms of both
  sides of the equation \eqref{eqn-10} and develop; we so get
\[ \sum \theta^\alpha \varphi^\beta
      \omega^\gamma \omega'^{\gamma'}
      \cdots \frac{1}{q^{1 + \varrho}}
   + {\textstyle\frac{1}{2}}
      \sum \theta^{2\alpha} \varphi^{2\beta}
      \omega^{2\gamma} \omega'^{2\gamma'}
      \cdots \frac{1}{q^{2 + 2\varrho}}
   + \cdots = \log L, \]
where the indices $\alpha, \beta, \gamma, \gamma',\ldots$ are with
respect to~$q$, as well as the sum. If we express the roots
$\theta, \varphi, \omega, \omega',\ldots$ in the way mentioned in
\S.~8 and put:
\[ \theta  = \Theta^{\mathfrak{a}},\qquad
   \varphi = \Phi^{\mathfrak{b}},\qquad
   \omega  = \Omega^{\mathfrak{c}},\qquad
   \omega' = \Omega'^{\mathfrak{c}'},\ldots,\]
then the general term on the left side becomes:
\[ \frac{1}{h} \sum \Theta^{h \alpha \mathfrak{a}}
                    \Phi^{h \beta \mathfrak{b}}
                    \Omega^{h \gamma \mathfrak{c}}
                    \Omega'^{h \gamma' \mathfrak{c}'}
      \cdots \frac{1}{q^{h + h\varrho}},\]
while after \eqref{eqn-11} we have to write for the right side:
\[ \log L_{\mathfrak{a}, \mathfrak{b}, \mathfrak{c}, \mathfrak{c}' \cdots} \]

Let now $m$ be an integer smaller than and coprime to $k$. If we multiply
both sides with:
\[ \Theta^{-\alpha_m \mathfrak{a}}
      \Phi^{-\beta_m \mathfrak{b}}
      \Omega^{-\gamma_m \mathfrak{c}}
      \Omega'^{-\gamma_m' \mathfrak{c}'}
      \cdots \]
and abbreviate the left side to the general term we get:
\begin{eqnarray*}
\cdots + \frac{1}{h} \sum \Theta^{(h \alpha - \alpha_m) \mathfrak{a}}
                    \Phi^{(h \beta - \beta_m) \mathfrak{b}}
                    \Omega^{(h \gamma - \gamma_m) \mathfrak{c}}
                    \Omega'^{(h \gamma' - \gamma_m') \mathfrak{c}'}
      \cdots \frac{1}{q^{h + h\varrho}} + \cdots
   \hspace{-216pt}\\
   &=& \Theta^{-\alpha_m \mathfrak{a}}
       \Phi^{-\beta_m \mathfrak{b}}
       \Omega^{-\gamma_m \mathfrak{c}}
       \Omega'^{-\gamma_m' \mathfrak{c}'}
       \cdots
   \log L_{\mathfrak{a}, \mathfrak{b}, \mathfrak{c}, \mathfrak{c}' \cdots}
\end{eqnarray*}
Summing now, to include all root combinations, from:
\[ \mathfrak{a}  = 0,\qquad
   \mathfrak{b}  = 0,\qquad
   \mathfrak{c}  = 0,\qquad
   \mathfrak{c}' = 0,\ldots \]
to:
\[ \mathfrak{a}  = 1,\quad
   \mathfrak{b}  = 2^{\lambda - 2} - 1,\quad
   \mathfrak{c}  = (p - 1) p^{\pi - 1} - 1,\quad
   \mathfrak{c}' = (p' - 1) p'^{\pi' - 1} - 1,\ldots,\]
the general term on the left hand becomes:
\[ \frac{1}{h} \sum W \frac{1}{q^{h + h\varrho}},\]
where the sum is over all primes~$q$ and $W$ means the product of
the sums taken over
$\mathfrak{a}, \mathfrak{b}, \mathfrak{c}, \mathfrak{c}',\ldots$
or respectively over:
\[ \sum \Theta^{(h \alpha - \alpha_m) \mathfrak{a}},\qquad
   \sum \Phi^{(h \beta - \beta_m) \mathfrak{b}},\qquad
   \sum \Omega^{(h \gamma - \gamma_m) \mathfrak{c}},\qquad
   \sum \Omega'^{(h \gamma' - \gamma_m') \mathfrak{c}'},\ldots.\]
We can now see from \S.~7 that the first of these sums is $2$ or $0$,
corresponding to if the congruence
$h \alpha - \alpha_m \equiv 0$ (mod.~$2$) or, equally, the congruence
$q^h \equiv m$ (mod.~$4$) holds or not; that the second is $2^{\lambda - 2}$
or $0$ corresponding to if the congruence
$h\beta - \beta_m \equiv 0$ (mod.~$2^{\lambda - 2}$) or, equally,
the congruence $q^h = \pm m$ (mod.~$2^\lambda$) holds or not;
that the third is $(p - 1) p^{\pi - 1}$ or $0$, corresponding to if
the congruence $h \gamma - \gamma_m \equiv 0$ (mod.~$(p-1)p^{\pi - 1}$)
or, equally, the congruence $q^h \equiv m$ (mod.~$p^\pi$) holds or not,
and so on; that therefore $W$ always disappears except when the congruence
$q^h \equiv m$ holds modulo each of the modules
$2^\lambda, p^\pi, p'^{\pi'},\ldots$, that is, when
$q^k \equiv m$ (mod.~$k$) holds, in which case $W = K$.
Our equation thus becomes:
\begin{eqnarray}
\label{eqn-15}
\sum \frac{1}{q^{1 + \varrho}}
      + {\textstyle\frac{1}{2}} \sum \frac{1}{q^{2 + 2\varrho}}
      + {\textstyle\frac{1}{3}} \sum \frac{1}{q^{3 + 3\varrho}}
      + \cdots \hspace{-144pt} \nonumber\\
   &=& \frac{1}{K} \sum \Theta^{-\alpha_m \mathfrak{a}}
       \Phi^{-\beta_m \mathfrak{b}}
       \Omega^{-\gamma_m \mathfrak{c}}
       \Omega'^{-\gamma_m' \mathfrak{c}'}
       \cdots
   \log L_{\mathfrak{a}, \mathfrak{b}, \mathfrak{c}, \mathfrak{c}' \cdots},
\end{eqnarray}
where the summation on the left is over all primes~$q$ the first, second,
third powers of which are contained in the form $\mu k + m$, while the
summation on the right is over
$\mathfrak{a}, \mathfrak{b}, \mathfrak{c}, \mathfrak{c}',\ldots$
and extends between the given limits. For $m= 1$, we get $\alpha_m = 0$,
$\beta_m = 0$, $\gamma_m = 0$, $\gamma_m' = 0$,\dots, and the right side
reduces to:
\[ \frac{1}{K} \sum
   \log L_{\mathfrak{a}, \mathfrak{b}, \mathfrak{c}, \mathfrak{c}' \cdots}.\]
The term of this sum corresponding to the $L$ of the first class,
$L_{0,0,0,0,\ldots}$, will, because of the expression \eqref{eqn-12}, 
contain $ \log \bigl( \frac{1}{\varrho} \bigr)$. Those terms corresponding
to different $L$ of the second class will, on condition of the desired
proof above, remain finite for infinitely small $\varrho$. Would the
limit for an arbitrary $L$ of the third class be zero the study of
the continuity, as in \S.~5, of the expression \eqref{eqn-13} for the
logarithm of this $L$, together with its $L$, would result in the term:
\[ -2 \log \left( \frac{1}{\varrho} \right) \]
which, combined with
$ \log \bigl( \frac{1}{\varrho} \bigr)$
in $\log L_{0,0,0,0,\ldots}$ would result in
$ - \log \bigl( \frac{1}{\varrho} \bigr)$
which would become $-\infty$ for infinitely small $\varrho$, while the
left side would consist of only positive terms. Therefore no $K$ of the
third class can have the limit zero, and it follows (excepting the missing
proof for the series of second class) that:
\[ \log L_{\mathfrak{a}, \mathfrak{b}, \mathfrak{c}, \mathfrak{c}' \cdots} \]
always approaches a finite limit for infinitely small $\varrho$, except
when simultaneously
$\mathfrak{a}  = 0$,
$\mathfrak{b}  = 0$,
$\mathfrak{c}  = 0$,
$\mathfrak{c}' = 0$,\dots
in which case the logarithm gets a value that is infinitely large.

Applying this result to the general equation \eqref{eqn-15} we see at once
that its right side becomes infinite for infinitely small $\varrho$,
namely by the term
$ \frac{1}{K} \log L_{0,0,0,0,\ldots}$
which grows beyond all limits, while all other remain finite.
Therefore also the left side must exceed any finite limit, from which
follows, as in \S.~6, that the series $\sum \frac{1}{q^{1 + \varrho}}$
has infinitely many terms or, in other words, that the number of
primes~$q$ of form $k\mu + m$, with $\mu$ an arbitrary integer and
$m$ a given number coprime to $k$, is infinite q.~e.~d.

\bigbreak

\centerline{\S.~11.}

\nobreak\bigskip

Regarding now the demonstration necessary for the completion of the
just developed proof, it reduces, according to the expression given
in \eqref{eqn-14} for the limit of an $L$ of second or third class, to
showing that for any of the root combinations of form:
\[ \pm 1,\, \pm 1,\, \pm 1,\, \pm 1,\,\ldots,\]
with the only exception of:
\[ +1,\, +1,\, +1,\, +1,\,\ldots,\]
the sum:
\begin{equation}
\label{eqn-16}
\sum (\pm 1)^\alpha (\pm 1)^\beta (\pm 1)^\gamma (\pm 1)^{\gamma'}
      \cdots \frac{1}{n},
\end{equation}
has a nonzero value, where $\alpha, \beta, \gamma, \gamma',\ldots$
means the system of indices for $n$, and where for $n$ are substituted
all positive integers not divisible by any of the primes $2, p, p', p'',\ldots$,
ordered by size. In the originally presented paper I proved this property
using indirect and quite complicated considerations. Later however I
convinced myself that the same object can be reached otherwise far shorter.
The principles from where we started can be applied to several other
problems, between which and the subject here treated one first would
guess to be no connection. For example, we could solve the very
interesting task to determine the number of different quadratic forms
which correspond to an arbitrary positive or negative determinant, and
we find that this number (which however is not the end result of this
investigation) can be expressed as product of two numbers, the first of
which is a very simple function of the determinant that has a finite
value for every determinant, while the other factor is expressed by
a series that coincides with the above \eqref{eqn-16}. From this then
follows immediately that the sum \eqref{eqn-16} never can be zero
since otherwise the number of quadratic forms for the respective
determinant would reduce to zero, while this number actually always
is $\geq 1$.

For this reason I will omit my earlier proof for the said property
of the series \eqref{eqn-16} here, and refer on the subject to the
mentioned investigation on the number of quadratic
forms\footnote{Preliminary notice on this can be found in \cite{D}
from which emerges the necessary theorem as corollary.}.

\end{document}